
\input amstex
\documentstyle{amsppt}

\magnification=\magstep1
\input epsf

\define \T {{\widetilde{T\Bbb{P}}{}^1}}
\define \dd#1{\frac{\partial}{\partial#1}}
\define \deq{\overset\text{def}\to=}
\define \one{\Bbb{I}}

\topmatter
\title A closed form for unitons \endtitle
\author Christopher Kumar Anand \endauthor
\affil McGill University and University of Warwick \endaffil
\address  Mathematics Research Centre,
  University of Warwick,
  Coventry CV4 7AL\endaddress
\email anand\@maths.warwick.ac.uk \endemail
\date 18 December 1995\enddate
\thanks Research supported by NSERC and FCAR scholarships and Professor
J.C.\ Hurtubise. \endthanks
\abstract
Unitons, i.e.\ harmonic spheres in a unitary
group, correspond to \lq uniton bundles\rq, i.e.\ holomorphic bundles
over the compactified tangent space to the complex line with certain
triviality and other properties.  In this paper, we use a monad
representation
similar to
Donaldson's representation of instanton bundles to obtain a simple
formula for
the unitons.  Using the monads,
we show that real triviality for
uniton bundles is automatic.  We interpret the uniton number as the `length'
 of a jumping line of the bundle, and identify the uniton bundles which
correspond to based maps into Grassmannians.  We also show that energy-$3$
unitons
are $1$-unitons, and give some examples.  \endabstract
\toc \widestnumber\head{19.}
\head 1. Introduction\page 2\endhead
\head 2. Prerequisites\page 4\endhead
\head 3. Uniton bundles and $\Bbb{P}^{2}$ monads\page 6\endhead
\head 4. Formula\page {10}\endhead
\head 5. Example of a $\operatorname{U}(3)$ uniton\page {13}\endhead
\head 6. Interpretation of Uniton Number\page {14}\endhead
\head 7. Real Triviality\page {16}\endhead
\head 8. Time Translation\page {17}\endhead
\head 9. Monads with $\gamma=0$\page {19}\endhead
\head 10. Example: $\operatorname{Harm}(\Bbb S^{2},
\operatorname{U}(2))$\page {20}\endhead
\head 11. Grassmannian Solutions\page {24}\endhead
\endtoc
\endtopmatter

\document
\NoBlackBoxes
\pageheight{9 true in} 
\pagewidth{6.5 true in}

\newpage
\head 1. Introduction\endhead
Harmonic maps between Riemannian manifolds $M$ and $N$
are critical values of an energy functional
$$\text{energy}(S:M\to N)=\frac12\int_M|dS|^2.$$
In the case of surfaces in a matrix group, with
the standard (left-invariant) metric, the energy takes the form
$$\text{energy}(S)=\frac12\int_{\Bbb R^2}
\left(|S^{-1}\dd{x} S|^2+|S^{-1}\dd{y} S|^2\right) dx\wedge
dy.\tag{1.1}$$
{\it Unitons} are harmonic maps  $S:{\Bbb S}^2\to \operatorname{U}(N)$.
Some authors call them multi-unitons.  Since the energy is
conformally invariant in the case of surfaces, it is natural to use
coordinates $x$ and $y$ on $\Bbb R^{2}$ (or $z\in \Bbb C$)
and derive the Euler-Lagrange or
uniton equations,
$$\dd{x} (S^{-1}\dd{x} S)+\dd{y} (S^{-1}\dd{y} S)=0.\tag{1.2}$$

{}From \cite{{SaUhl}, Theorem 3.6}, we know that
harmonic maps from $\Bbb R^2\to \operatorname{U}(N)$ extend to
${\Bbb S}^2$ iff they
have
finite
energy, and that such maps are always smooth.
So working in terms of coordinates
$x$ and $y$ on $\Bbb R^2$ or $z\in\Bbb C$ poses no real limitation.

\subhead {1.3} Based Unitons \endsubhead
Unitons are determined by the pullback of the Maurer-Cartan form on
$\operatorname{U}(N)$,
$$A\deq\frac12S^{-1}d S=A_xdx+A_ydy=A_zdz+A_{\bar{z}} d\bar{z}
\tag{1.4}$$
and a choice of a basepoint, $S(\infty)\in \operatorname{U}(N)$,
as we can see by
thinking
of $d+2A$ as a flat connection and $S$ as a gauge transformation.
(See \cite{An1}.)

We are concerned with the based maps
$$\operatorname{Harm}^{*}_{k}(\Bbb{S}^{2},
\operatorname{U}(N))\deq\left\{S\in
\operatorname{Harm}(\Bbb{S}^{2}, \operatorname{U}(N)):S(\infty)=1,
\text{energy}(S)=k\right\}.$$

In \cite{Uhl}, Uhlenbeck showed that all unitons could be constructed
from simpler unitons by `adding a uniton'.  This construction was
investigated, from different
perspectives, by Wood, Valli, Guest, Ohnita and Segal.
We approach the question of constructing unitons and
investigating their moduli via algebraic integration, using a twistor
construction of Hitchin and Ward (\cite{Hi},\cite{Wa}).  We
proved in \cite{An1} that the based unitons,
$\operatorname{Harm}^{*}(\Bbb{S}^{2},
\operatorname{U}(N))$,
 are isomorphic to uniton bundles,
with energy corresponding to the bundles' second
Chern class.
Given transition matrices for the uniton bundle, we
 showed how to construct solutions. This equivalence is explained in \S 2.

In this paper we apply Horrocks' monad construction to the uniton
bundles, i.e.\ we show that uniton bundles are representable as the
cohomology of a short sequence (of homogeneous bundles on $\Bbb{P}^{2}$
in this case).   After interpreting the reality and triviality properties
in \S 3, we derive the main result of this paper,
\proclaim{Theorem A} Based, rank-$N$ unitons of energy $8\pi k$ are all of
 the form
 $$S=\Bbb I +a\alpha_{2}^{-1}(\alpha_{1}-2(x+iy\alpha_{2}))^{-1}b$$
 for some choice of $N\times k$,
 $k\times N$, $k\times k$ and $k\times k$ matrices,
 $a$, $b$, $\alpha_{1}$, $\alpha_{2}$.  (Multiplication is matrix
 multiplication.)
\endproclaim \noindent
We do this
by specialising the `monodromy' interpretation of Uhlenbeck's
 extended solution (see \S2.11)
to bundles given by monads.
In \S{5} we explain how to obtain the monad
representation for the example of a $2$-uniton in \cite{An1}. In \S 6,
we show that the uniton number is also the length of the jumping line which
determines the uniton bundle.  We prove in the next section that real
triviality
is implied by the other bundle properties.  This is a geometric fact whose
proof depends on the monad representation.  We complete the GIT picture
in \S 8 by representing time invariance with its infinitesimal generator.
As
a corollary, we give a simple finite-gap result.  The next two sections work
out the $1$-uniton case in general.  Finally, in \S 11, we explain how to
identify those solutions which factor through the Cartan embedding of a
Grassmannian.  Precisely, we have

\proclaim{Theorem B} \roster
\item The space of based unitons $\operatorname{Harm}^{*}_{k}
(\Bbb{S}^{2},\operatorname{U}(N))$ is isomorphic to the set of monad
data
$$\gather \gamma, \alpha'_{1},\delta \in\operatorname{gl}(k/2),\quad \gamma
\text{ nilpotent}\\
a'\in\operatorname{M}_{N,k/2},\quad
b'\in\operatorname{M} _{k/2,N}\endgather$$
satisfying
$$\gather \operatorname{rank} \pmatrix\gamma\\ \alpha'_{1}+z\\ a'
\endpmatrix
=\operatorname{rank}\pmatrix\gamma&\alpha'_{1}+z&b'
\endpmatrix=k/2\qquad\forall z\in\Bbb{C}\tag{\text{nondegeneracy}}
\\
[\alpha'_{1},2\gamma ]+b'a'=0\tag{\text{monad equation}}
\endgather$$
$$\matrix [\delta,\gamma]=0 & a'\delta=0\\
[\delta,\alpha'_{1}]=\gamma & \delta b'=0\endmatrix\tag{\text{time
invariance}}$$
quotiented by the action
of $g\in\operatorname{Gl}(k/2)$
$$\gathered \gamma\mapsto g\gamma g^{-1}\quad
\alpha'_{1}\mapsto g\alpha'_{1} g^{-1}
\quad \delta\mapsto g\delta g^{-1}\\
a'\mapsto a' g^{-1}\quad b'\mapsto gb'.\endgathered\tag{group action}$$
\item These data determine the uniton bundle over a hemisphere.  Reality
determines it over the other hemisphere, giving monad data as in
Theorem A as follows:
$$\matrix \alpha_{1}=2\pmatrix\alpha_{1}^{\prime *}& \phi_{1}\\
\phi_{2}&\alpha'_{1}
\endpmatrix&
\alpha_{2}=\pmatrix -\one-2\gamma^{*}&\\&\one+2\gamma
\endpmatrix\\
a=2\pmatrix b^{\prime *}&a'
\endpmatrix&
b=2\pmatrix a^{\prime *}\\b'
\endpmatrix\endmatrix\tag{\text{reality}}$$
where $\phi_{1}$ and $\phi_{2}$ are functions of $\gamma$, $a'$ and
$b'$ determined by the big monad equation
$[\alpha_{1},\alpha_{2}]+ba=0$.
\item The uniton number (see \S6)
   is the smallest $n\in\Bbb{Z}$ such that $\gamma^n=0$.
\item The uniton is in the image of $\operatorname{Harm}^{*}_{k} (\Bbb S^2,
\operatorname{Gr}(\Bbb C^N))\hookrightarrow \operatorname{Harm}^{*}_{k}
(\Bbb S^2,\operatorname{U}(N))$ iff there is a $G\in \operatorname{Gl}(k/2)$,
and an $F\in\operatorname{U}(N)$ such that $F^2=\Bbb I$ and
$$\matrix G\alpha_{2}^{-1}\alpha_{1}G^{-1}=\alpha_{1} &
Fa\alpha_{2}^{-1}G^{-1}=a\\
G\alpha_{2}G^{-1}=\alpha_{2}&
-G\alpha_{2}^{-2}bF^{-1}=b.
\endmatrix$$
\endroster
\endproclaim

Clearly Theorem B can be used to investigate the topology of the moduli
space of unitons (see \cite{An2} for some results in low dimensions),
but the details of such calculations properly belong
in a separate paper.  In this paper, we restrict ourselves to showing that
the space of uniton bundles has high codimension as a subspace of the space of
framed jumping lines: generic
nilpotent  $\gamma$ are not allowed.  This implies that 2-unitons must
have normalised energy 4 or more.  It would be interesting to formulate
and prove a more general finite-gap result of this type.

\head 2. Prerequisites\endhead
\subheading{{2.1} Uniton Bundles}
The uniton bundles are bundles on $\widetilde{T\Bbb{P}}{}^{1} \deq
\Bbb{P}{}(\Cal O{}\oplus\Cal O{}(2))$,
the fibrewise compactification of the tangent bundle $T\Bbb{P}{}^{1}$ of
the
{\it complex}
projective line.

Let $(\lambda,\eta)$  and
$(\hat{\lambda}=1/\lambda,\hat{\eta}=\eta/\lambda^2)$ be coordinates
on $T\Bbb{P}^1\cong
\Cal{O}_{\Bbb P^{1}} (2)$,
where $\lambda$ is the usual coordinate on $\Bbb{P}^1$ and $\eta$
is the coordinate associated to $d/d\lambda$.
Meromorphic sections ($s$) of $T\Bbb{P}^1$ give all the holomorphic
sections of $\widetilde{T\Bbb{P}}{}^{1}$ ($[s,1]$ in projective
coordinates on
$\widetilde{T\Bbb{P}}{}^{1})$, save one.  We fix notation for the lines on
$\T$:
$$\aligned P_\lambda &=\pi^{-1}(\lambda\in\Bbb{P}^1)=\text{a pfibre
(silent p)}\\ {G}_0&=\{(\lambda ,[0,1])\}=\text{(graph of) zero section of
$T\Bbb{P}^1$}\\
{G}_\infty&=\{(\lambda ,[1,0])\}=\text{infinity section of
$\widetilde{T\Bbb{P}}{}^{1}$}\\
{G}_{\eta=s}&=\left\{(\lambda
,[s(\lambda),1])\right\}.\endaligned\tag{2.2}$$
If $y=(a,b,c)\in\Bbb C^3$, we will also write ${G}_y$ for
${G}_{\eta=\frac12(a-2b\lambda -c\lambda^2)}$.

To encode unitarity, we need the real structure
$$\sigma^*(\lambda,\eta)=(1/\bar{\lambda} ,-\bar{\lambda}
^{-2}\bar{\eta})\tag{2.3}$$
which acts by
$$\sigma^*(a,b,c)=(\bar c,-\bar
b,\bar a)$$
on $\Bbb C^3\cong H^{0}(\Bbb{P}^{1},\Cal{O}(2))$, the space of finite
sections.
We similarly define time translation
$$\align\delta_{t}:\ & (\lambda,\eta)\mapsto (\lambda,\eta -2t\lambda)\\ &
(a,b,c)\mapsto (a,b+t,c).\tag{2.4}\endalign$$

\proclaim{Definition {2.5}} A rank $N$, or $\operatorname{U}(N)$, {\rm uniton
bundle},
$\Cal{V}$,
is a holomorphic rank $N$ bundle on $\widetilde{T\Bbb{P}}{}^{1}$ which is
a) trivial when restricted
to the following curves
in $\widetilde{T\Bbb{P}}{}^{1}$
\roster
\item the section at infinity
\item nonpolar fibres (i.e.\ fibres above $\lambda\in\Bbb
C^{*}\subset\Bbb{P}^1$)\item real sections of $T\Bbb{P}^1$
(sections invariant under $\sigma$)
\endroster
b) is equipped with bundle lifts
$$\CD \Cal{V}@>\tilde\delta_t>>\Cal{V}\\@VVV
@VVV\\\widetilde{T\Bbb{P}}{}^{1}@>\delta_t>>\widetilde{T\Bbb{P}}{}^{1}
\endCD\quad\text{and}\quad
\CD \Cal{V}@>\tilde\sigma>>\Cal{V}^*\\@VVV @VVV\\
\widetilde{T\Bbb{P}}{}^{1}
@>\sigma>>\widetilde{T\Bbb{P}}{}^{1}\endCD$$
\roster
\item $\tilde\delta_t$ a one-parameter family of holomorphic
transformations fixing $\Cal{V}$ above the section at infinity,
lifting $\delta_t$,  and
\item $\tilde\sigma$ a norm-preserving, antiholomorphic lift of
 $\sigma$ such that the induced
hermitian metric on $\Cal{V}$ restricted to a fixed point of $\sigma$
is positive definite; equivalently, such that the induced lift to  the
principal bundle of frames  acts on fibres of fixed points of $\sigma$ by
$X\mapsto X^{*-1}$.\endroster
and c) has a framing,
$\phi\in H^0(P_{-1},\text{frames}(\Cal{V}))$, of the bundle
$\Cal{V}$ restricted
to the fibre $P_{-1}=\{\lambda=-1\}\subset \widetilde{T\Bbb{P}}{}^{1}$
such that $\tilde\sigma(\phi)=\phi$.
\endproclaim

\subhead {2.6} Extended Solutions\endsubhead
Uhlenbeck's extended solutions $E_\lambda$
(actually
first employed in \cite{Po}), encode the unitons as follows
\proclaim{Theorem {2.7} \cite{{Uhl}, 2.1}}  Let $\Omega\subset {\Bbb
S}^2$ be a simply-connected
neighbourhood and
$A:\Omega\to T^*(\Omega)\otimes \operatorname{u}(N)$.
Then $2A=S^{-1}d S$, with $S$ harmonic iff the curvature of the
connection $$\Cal{D}_\lambda =(\dd{\bar{z}}+(1+\lambda)A_{\bar{z}},
\dd{z} +(1+\lambda^{-1})A_z)\tag{2.8}$$
vanishes for all $\lambda\in\Bbb C^*$.\endproclaim
\proclaim{Theorem {2.9} \cite{{Uhl}, 2.2}} If $S$ is harmonic and
$S(\infty)=\Bbb I$, then there exists a unique holomorphic family of
covariant constant
frames $E_\lambda :\Bbb{S}^2\to
\operatorname{U}(N)$ for the connection $\Cal{D}_\lambda$ for each
$\lambda\in\Bbb
C^*$ with
\roster \item $E_{-1}=\Bbb I$, \item $E_1=S$,
\item
$E_\lambda(\infty)=\Bbb I$.\endroster
  Moreover, $E_\lambda$ is analytic and holomorphic in
$\lambda\in\Bbb C^{*}$.\endproclaim
\proclaim{Theorem {2.10} \cite{{Uhl}, 2.3}}
Suppose $E:\Bbb C^{*}\times\Omega\to \operatorname{Gl}(N)$ is
analytic and holomorphic in the first variable, $E_{-1}\equiv\Bbb I$, and
the expressions
$$\frac{E_\lambda^{-1}\bar{\partial} E_\lambda}{1+\lambda},
\quad\frac{E_\lambda^{-1}\partial
E_\lambda}{1+\lambda^{-1}}$$
are constant in $\lambda$, then $S=E_1$ is harmonic.\endproclaim

The key to understanding the `monodromy' interpretation of $E_{\lambda}$
which we require from \cite{An1} is to think of $E_\lambda$ as a singular
change of frame between the usual (constant) frame of $\Bbb C^N$ and the
trivialising frame for $\Cal D_\lambda$ which agrees with first frame at
$z=\infty$.

\subhead {2.11} Monodromy construction \endsubhead
The uniton bundle $\Cal{V} \to \widetilde{T\Bbb{P}}{}^{1}$ is constructed
generically as the kernel of a differential operator on $\Bbb R^{3}
\times \Bbb S^{2}$ (\cite{An1}).  The `connection'
$\Cal{D}_{\lambda}$  is the projection of this operator from
$\Bbb R^{3}\times\Bbb S^{2}$ to $\Bbb R^{2}\times\Bbb C^{*}$.  Pulling
$E_{\lambda}$ back to $\Bbb R^{3}\times\Bbb C^{*}$ gives a solution to the
defining operator on $(\Bbb S^2 \times \Bbb R)\times\Bbb C^*$ which gives
a trivialisation
of $\Cal{V}$ restricted to the open set $\{\lambda\in\Bbb C^{*}\}$.
Intrinsically, it is the trivialisation on each fibre $P_{\lambda}$,
$\lambda\in\Bbb C^{*}$, which agrees with a trivialisation of $\Cal{V}
|_{G_{\infty}}$.  We require this trivialisation to agree with the framing
$\phi$ at $P_{-1}$.  By construction, the usual frame above a point $y\in\Bbb
R^3$
also lifts to give a solution of the operator, this time giving a frame over
the section $G_y$.  This frame is also
uniquely defined by compatibility with the
framing.  The change of frame $E_\lambda$ can be computed
by composing the cycle of isomorphisms
$$\CD
\Cal{V}_{\lambda,\infty} @<\text{eval}<< H^0({G}_\infty,\Cal{V})
@>\text{eval}>>
  \Cal{V}_{-1,\infty}\\
@A\text{eval}AA @. @AA\text{eval}A\\
H^0(P_\lambda,\Cal{V}) @. @. H^0(P_{-1},\Cal{V})\\
@V\text{eval}VV @. @VV\text{eval}V\\
\Cal{V}_{(\lambda,z/2-t\lambda-\lambda^2\bar{z}/2)} @<\text{eval}<<
H^0({G}_{(z,\bar{z},t)},\Cal{V})
  @>\text{eval}>> \Cal{V}_{(-1,z/2-t\lambda-\bar{z}/2)}\endCD\tag{2.12}$$
counter clockwise.  The existence of the bundle isomorphism
$\tilde{\delta}_{t}$ (time translation) ensures that the result doesn't
depend on $t$.  Finiteness, i.e.\ extension to $\Bbb S^{2}$, follows from
the compactness of $\T$.

\head 3. Uniton Bundles and $\Bbb{P}^{2}$ monads\endhead

Uniton bundles are bundles over $\T$, the $\Bbb{P}^{1}$ bundle over
$\Bbb{P}^{1}$ with a double twist.  In \cite{An1}, the geometry of $\T$
allowed us to construct the uniton bundle from the uniton and the extra
structure (reality, time invariance) fits naturally with this geometry.
In this case, however, what appears natural is not entirely optimal.  By
transporting the uniton bundles to $\Bbb{P}^{2}$ via a birational
equivalence, constructing and manipulating monads becomes much easier.
To start with, the existence of monad representations is known.

Every operation we will make on $\Bbb{P}^{2}$ has its analogue
on $\T$, but monads on projective spaces have two big advantages which we
will exploit in interpreting the reality and triviality
properties of the bundle.  Namely,
\roster
\item $\Bbb P^{2}$ monads are self-dual, i.e.\ the transposed monad is a monad
of the same form representing the dual bundle, and
\item we know when the bundle is trivial on hyperplanes.  In the usual
notation (described below), $\Cal V\to \Bbb P^{2}$ is trivial on a hyperplane
$L=\overline{p_{1}p_{2}}\iff\det (K_{p_{2}}\circ J_{p_{1}})\ne 0$, and any
choice of spanning representatives of $\ker K_{p_{2}}/\operatorname{im}
J_{p_{1}}$
frames
$\Cal V|_{L}$ canonically, in particular a basis for
$\ker K_{p_{1}}\cap\ker K_{p_{2}}$ does.
(See \cite{OSS}, \cite{Do}, \cite{Hu}.)
\endroster

\subhead 3.1 The birational equivalence\endsubhead
If $X$, $Y$ and $W$
are homogeneous coordinates on $\Bbb P^{2}$, and $\lambda$ and $\eta$
base and fibre coordinates on $\T$ (see \thetag{2.1}),
$$\{ X=\lambda Y,\ (X+Y)W=\eta Y^{2}\}\subset \T\times\Bbb P^{2}\tag{3.2}$$
is the graph of the birational equivalence which comes from
\roster \item blowing up the point $(\lambda =-1,\eta=0)$,
\item blowing down $\widetilde{P}_{-1}$ (the
proper transform of the fibre $\{\lambda=-1\}$), and
\item blowing down the image of $G_{\infty}$.  \endroster
The birational equivalences we will need are given
diagramatically in {figure 1}.

\midinsert
\hfill\epsfxsize=.9\hsize\epsffile{blowups.epsf}\hfill
\botcaption{Figure 1}
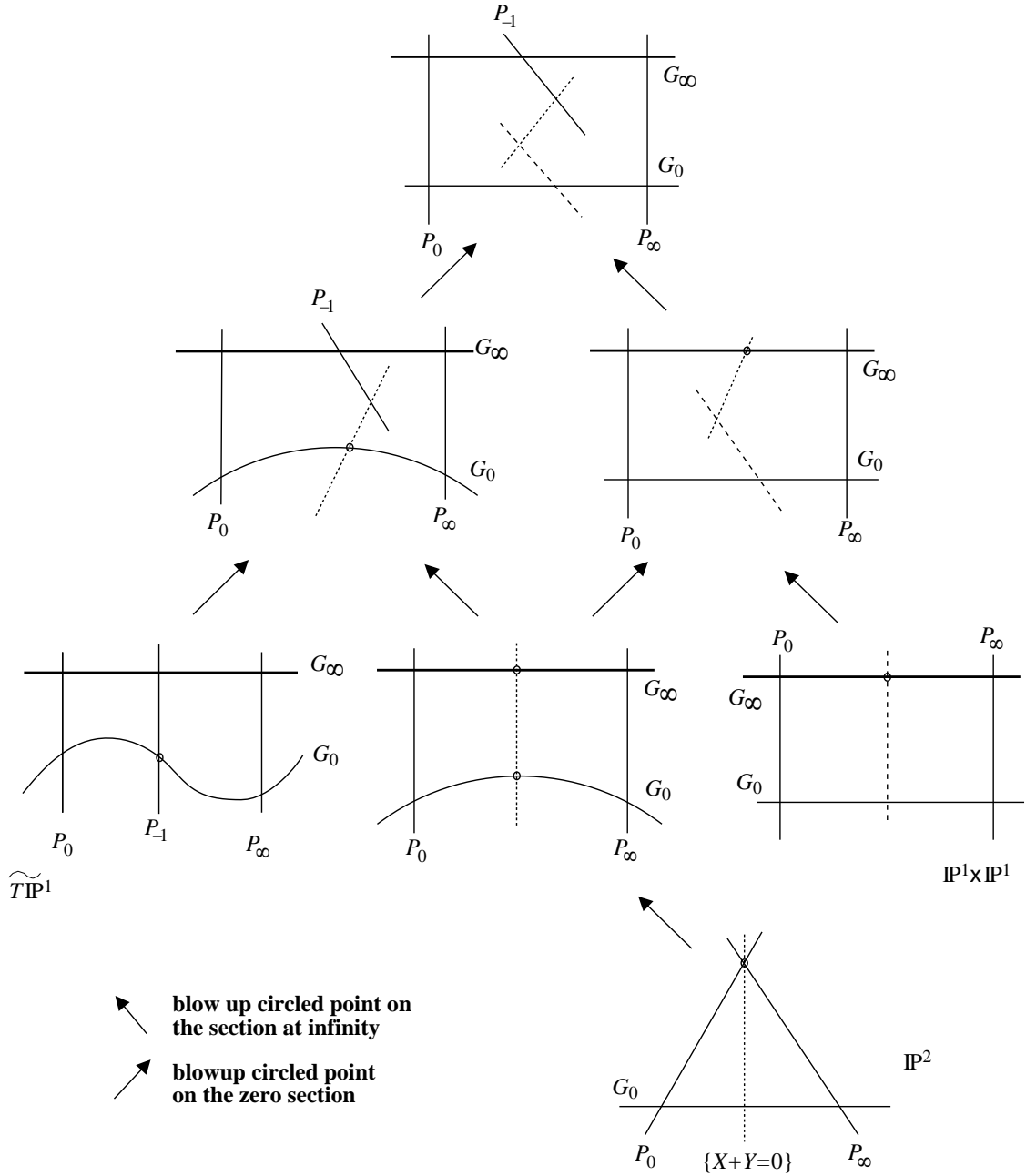
Birational equivalences of surfaces can be realised by sequences
of blowings up and blowings down.   Monad representations for
generically trivial bundles
 exist for all rational surfaces and monad maps equating
these representations are defined on a `common ancestor'
of the two surfaces.
For example in the explicit monad construction of \S5 we need to work on
the topmost surface.
\endcaption
\endinsert

Under the birational equivalence \thetag{3.2}, the ruling $\{
P_{\lambda}:\lambda\in\Bbb
P^{1}\}$ of $\T$
is mapped to the pencil of lines $\{X=\lambda Y\}$ on $\Bbb P^{2}$,
except the fibre $P_{-1}$, which is mapped to a point.  On $\Bbb{P}^{2}$, the
line $X+Y=0$ is the exceptional divisor, and $G_{\infty}$ becomes the point
$[0,0,1]\in\Bbb{P}^2$.  Since push forward gives an isomorphism (\cite {{An1},
Lemma 7{}.{}2},\cite {At})
between bundles on $\T$ and
bundles on $\Bbb P^{2}$, trivial on $P_{-1}$ and $\{X+Y\}$ respectively,
we will use the same letter for a bundle and its pullback.

Assume now that $\Cal V$ is a uniton bundle.  Since ${\Cal V}$ is trivial
on generic lines, $\Cal V$
admits a monad representation (\cite {{OSS}, example 3, p249}),
$\Cal V = \ker K/\operatorname{im} J$,
where
$$0\to\Cal O_{\Bbb P^{2}}(-1)^{k} @>J>> {\Cal O_{\Bbb P^{2}}}^{2k+N}
@>K>> \Cal O_{\Bbb P^{2}}(1)^{k}\to 0\tag{3.3}$$
is a complex of linear maps, i.e.\ degree one homogeneous polynomials
such that
$K\circ J=0$, on each fibre $J$ is injective and $K$ is surjective,
and $k=c_{2}\Cal{V}$.

Since the birational equivalence blows down the fibre $P_{-1}$ on which
$\delta_{t}$ acts nontrivially and replaces it with an exceptional fibre on
which $\delta_{t}$ must act discontinuously, one must either work with $\T$
monads or use another birational transformation, which we do in \S{8}.

Clearly the monad representation \thetag{3.3} is not unique;
$\operatorname{Gl}(k)\times \operatorname{Gl}(2k+N)\times
\operatorname{Gl}(k)$ acts on the vector spaces linearly inducing monad
equivalences.  Given that
$\Cal V|_{\{X+Y\}}$ is trivial, we can assume that the monad has the
block form
$$\matrix J_{W}=\pmatrix \Bbb I\\0\\0\endpmatrix &
J_{X}-J_{Y}=\pmatrix 0\\ \Bbb I\\0\endpmatrix &
J_{X}+J_{Y}=\pmatrix \alpha_{1}\\ \alpha_{2} \\a\endpmatrix \\
\qquad&\qquad&\qquad\\
K_{W}=\pmatrix 0 & \Bbb I & 0 \endpmatrix &
K_{X}-K_{Y}=\pmatrix -\Bbb I & 0 & 0 \endpmatrix &
K_{X}+K_{Y}=\pmatrix -\alpha_{2} & \alpha_{1} & b \endpmatrix
\endmatrix\tag{3.4}$$
where the $\alpha_{i}$, $a$ and $b$ are $k\times k$, $N\times k$ and $k\times
N$
matrices, respectively.  This form is stabilised by an action of
$\operatorname{Gl}(k) \times \operatorname{Gl}(N)$.  The
$\operatorname{Gl}(N)$ action corresponds to changes of frame.  Since uniton
bundles are framed it does not act.

The monad equation $K\circ J=0$ is
$$[\alpha_{1},\alpha_{2}]+ba=0\tag{3.5}$$
and nondegeneracy says that
$$\matrix \pmatrix \alpha_{1} + X\\ \alpha_{2}+Y \\ a \endpmatrix
& \text{and} & \pmatrix -\alpha_{2}-Y & \alpha_{1}+X & b \endpmatrix
\endmatrix$$
are respectively injective and surjective for all $X,Y$.  Since
$$K_{W}\circ (\lambda J_{X}+J_{Y})=\frac{1}{2}((1+\lambda)\Bbb I
-(1-\lambda)\alpha_{2}),$$
$\Cal V|_{\{X=\lambda Y\}}$ is holomorphically nontrivial iff $
\frac{1+\lambda }{1-\lambda}$ is an eigenvalue of $\alpha_{2}$.

Since only $\Cal V|_{\{X=0\}}$ and $\Cal V|_{\{Y=0\}}$ are nontrivial,
$$\alpha_{2}=\pmatrix -\Bbb I -2\gamma^{\prime} & \\ & \one+2\gamma
\endpmatrix\tag{3.6}$$
for some nilpotent matrices $\gamma,\gamma^{\prime}$.  (The twos will be
convenient later.)  If we correspondingly decompose
$$\alpha_{1}=\pmatrix \alpha_{1,11} & \alpha_{1,12}\\ \alpha_{1,21} &
\alpha_{1,22} \endpmatrix$$
the condition \thetag{3.5} implies that $\alpha_{1,12}$ and
$\alpha_{1,21}$ are determined by $\gamma$, $\gamma'$, $a$ and $b$.

\subhead {3.7} Reality \endsubhead
In interpreting reality, we are motivated by the geometry of the uniton
bundle.  The bundle $\Cal V$ is trivial on fibres other than $P_{0}$ and
$P_{\infty}$ and is specified by these jumping lines,
i.e.\ by holomorphic bundles
on neighbourhoods of $P_{0}$ and
$P_{\infty}$ with framings along $G_\infty$ which tell how to glue them
into the trivial bundle on $\T\setminus (P_0\cup P_\infty)$.

Since
the real structure fixes the framing of $\Cal V$ and exchanges $P_{0}$
and $P_{\infty}$, we know that $\Cal V|_{P_{0}}$ and $\Cal
V|_{P_{\infty}}$ are `conjugate'.  We can separate the bundle
into two parts by gluing neighbourhoods of $P_{0}$ and
$P_{\infty}$ into separate trivial bundles, obtaining $\Cal
V_{\text{north}}$ and $\Cal V_{\text{south}}$ say.  Then
$$\Cal V_{\text{north}}\cong \overline{\sigma^{*}{\Cal
V_{\text{south}}}^{\text{dual}}},$$
and the corresponding monads $(J^{\prime},K^{\prime}), (J^{\prime\prime},
K^{\prime\prime})$ will satisfy
$J^{\prime\prime}=\overline{\sigma^{*}K^{\prime t}}$, $K^{\prime\prime}=
\overline{\sigma^{*}J^{\prime t}}$.  Since
$$\align X & \mapsto\overline{Y}\\ \sigma^{*}: Y & \mapsto\overline{X}\\ W &
\mapsto
-\overline{W}\endalign$$
we obtain (after putting the monad back into normal form)
$$\align \alpha^{\prime\prime}_{1} & =\alpha^{\prime *}_{1}\\
\alpha^{\prime\prime}_{2} & =-\alpha^{\prime *}_{2},\ \text{i.e.\ }
   \quad\gamma=\gamma^{\prime *}\\
a^{\prime\prime} & =b^{\prime *}\\
b^{\prime\prime} & =a^{\prime *}, \endalign$$
in particular $\alpha^{\prime}_{i}$ and $\alpha_{i}^{\prime\prime}$ have the
same size.

Finally, the monad equation determines $\alpha_{1,12}$ and $\alpha_{1,21}$
to be
$$\align \alpha_{2,12} & =-\frac{1}{2}\sum\gamma^{*(i-1)}(\Bbb I
+\gamma^{*})^{-i}a^{\prime *}a^{\prime}(\Bbb I
+\gamma)^{-i}\gamma^{i-1}\\
\alpha_{2,21} & =\frac{1}{2}\sum\gamma^{i-1}(\Bbb I
+\gamma)^{-i}b^{\prime}b^{\prime *}(\Bbb I +\gamma^{*})^{-i}\gamma^{*(i-
1)}.\tag{3.8}\endalign$$
Summing up, a (real) uniton bundle is represented by a $\Bbb P^{2}$ monad
$$\align \alpha_{2} & =\pmatrix -\Bbb I-2\gamma^{*} & \\ & \Bbb I
+2\gamma \endpmatrix \\
\alpha_{1} & =\pmatrix \alpha_{1}^{\prime *} & \sum\cdots \\ \sum\cdots &
\alpha_{1}^{\prime} \endpmatrix\\
b & =\pmatrix a^{\prime *}\\b^{\prime}
\endpmatrix \\
a & =\pmatrix b^{\prime *} & a^{\prime}\endpmatrix .\endalign$$

\subhead {3.9} Real Triviality \endsubhead
Real sections of $\T$ are smooth, linearly-equivalent curves of
self-intersection two.  By blowing up and down twice, some of these curves
become singular, which seems unpleasant, but look again at the singular
sections.

Sections $\eta = a+b\lambda +c$ of $\T$ are parametrised by $\Bbb
C^{3}$.  The sections $G_{(z,t,\bar{z})}=\{\eta
=\frac{1}{2}(z-2t\lambda-\bar{z}\lambda^{2})\}$ for $t,z\in\Bbb{C}$ are
time translates of real sections.
Uniton bundles are required to be trivial on real sections.  Since they are
 invariant under time translation,  it is enough to
know that for each $z\in \Bbb C$,
 $\Cal V$ is trivial on $G_{(z,t,\bar{z})}$ for some $t\in\Bbb
C$.  The sections with $t=(\bar{z}-z)/2$ all
contain the point $(\lambda =-1,\eta =0)$.  So the proper transforms
of these sections are singular (the union of two curves).  Since they
still have self-intersection two, and are effective, they are unions of
two hyperplanes $(\bar{z}X-zY+2W)(X+Y)=0$, one of which is the exceptional
fibre.  So $\Cal V|_{G_{(z,t,\bar{z})}}$ is trivial iff $\Cal V$ is trivial
on the hyperplane $\bar{z}X-zY+2W=0$.  So $\Cal V$ is trivial on all real
sections iff
$$\gather0\ne\det \left(\frac{1}{2}(z-\bar{z})K_{W}+K_{X}+K_{Y}\right)
\left(\frac{1}{2}(-z-\bar{z})J_{W}+J_{X}-J_{Y}\right)\\
=\det \left(\frac{1}{2}(z-\bar{z})\alpha_{2}+\frac{1}{2}(z+\bar{z})\Bbb I
-\alpha_{1}\right)\tag{3.10}\endgather$$
for all $z\in\Bbb C$.  In \S{7} we will show that real triviality is implied
by reality, but it will be helpful to have the expression \thetag{3.10}
in the next
section.

\head 4. Formula\endhead
\subhead {4.1} Formula\endsubhead
Now we are in a position to use the monodromy construction of the uniton
\thetag {2.12},
to get a formula for the uniton in terms of the monad data $\alpha_{1}$,
$\alpha_{2}$, $a$, $b$.  It suffices to parametrise the sections of $\Cal
V|_{\{X=\lambda Y\}}$ and of $\Cal V|_{\{\bar{z}X-zY+2W=0\}}$ because the
section at infinity is collapsed to a point, and both the bundle and
corresponding uniton are
time-translation invariant.

Sections of $\Cal V|_{\{\bar{z}X-zY+2W=0\}}$ are parametrised by
$$\ker K_{p_{1}}\cap\ker K_{p_{2}}$$
where $\overline{p_{1}p_{2}}=\{\bar{z}X-zY+2W=0\}$.  Taking
$p_{1}=[2,2,z-\bar{z}]$ and $p_{2}=[2,-2,-z-\bar{z}]$,
$$H^{0}(\Cal V|_{\overline{p_{1}p_{2}}})\cong\ker\pmatrix K_{p_2}\\
K_{p_1}\endpmatrix =\ker \pmatrix -\Bbb I &
\frac{1}{2}(z-\bar{z})\Bbb I & 0 \\ -\alpha_{2} & \alpha_{1}-
\frac{1}{2}(z+\bar{z})\Bbb I & b
\endpmatrix.$$

Real triviality \thetag{3.10} says
$$\det\pmatrix -\Bbb I & iy\Bbb I \\ -\alpha_{2} & \alpha_{1}-x
\endpmatrix = (-1)^{k/2}\det (\alpha_{1}-(x+iy\alpha_{2}))\ne0$$
which implies that the kernel is
$$\left\{\pmatrix -iy(\alpha_{1}-(x+iy\alpha_{2}))^{-1}bs \\
-(\alpha_{1}-(x+iy\alpha_{2}))^{-1}bs \\ s \endpmatrix
\in\Bbb C^{2k+N}:s\in\Bbb C^{N}\right\}.\tag{4.2}$$

Sections of $\Cal V|_{\{X=\lambda Y\}}$ are similarly given by
$$\pmatrix (\alpha_{2}+\frac{\lambda+1}{\lambda-1}\Bbb
I)^{-1}bs\\0\\s\endpmatrix\in\Bbb C^{2k+N}\tag{4.3}$$
as $s$ varies in $\Bbb C^{N}$.

The evaluation map $H^0(\{X=\lambda Y\},\Cal V)\to\Cal V_p$ is given by
\thetag{4.3}$\mapsto$\thetag{4.3}$+\operatorname{im}J_p$.  To compute
$$H^0(\{\bar{z}X-zY+2W=0\},\Cal V)\to\Cal V_p \to H^0(\{X+\lambda Y\},\Cal
V)$$
we have to take the expression \thetag{4.2} and add an element of
$\operatorname{im}J_p$ to obtain a
representative of the form \thetag{4.3}.

At the point of intersection $p=[-\lambda,-1,\bar{z}\lambda-z]$
$${\ssize\frac{2}{1-\lambda}}J_p=
\pmatrix \alpha_{1}+\frac{z-\lambda\bar{z}}{\lambda -1} \\
\alpha_{2}+\frac{\lambda +1}{\lambda -1} \\ a \endpmatrix.$$
When $\lambda\in\Bbb C^{*}$, we can translate \thetag{4.2} into the form
of \thetag{4.3} by adding
$${\ssize\frac{2}{1-\lambda}}J_p\left(\pmatrix\alpha_{2}+
\frac{\lambda+1}{\lambda-1}\endpmatrix^{-1}
\pmatrix\alpha_{1}-(x+iy\alpha_{2})\endpmatrix^{-1}bs\right)$$
giving
$$\pmatrix
\left(-iy+(\alpha_{1}+\frac{z-\lambda\bar{z}}{\lambda-1})\right)
(\alpha_{2}+\frac{\lambda+1}{\lambda-1})^{-1}(\alpha_{1}
-(x+iy\alpha_{2}))^{-1}bs \\
0 \\
\left(\Bbb I +a(\alpha_{2}+\frac{\lambda+1}{\lambda-1})^{-1}(\alpha_{1}
-(x+iy\alpha_{2}))^{-1}b\right)s \endpmatrix.$$
Note that the resulting map $\Bbb C^{N}\to\Bbb C^{N}$ is linear.  Since
the parametrisations of $H^{0}(\{X=\lambda Y\}, \Cal V)$ agree at their
point of intersection $[0,0,1]$, we can compute the monodromy
construction of the extended solution as a composition of two linear
maps---a matrix:
$$\aligned E_{\lambda}=&\left(\one +a\left(\alpha_{2}+
  \matrix{\frac{\lambda +1}{\lambda -1}}\endmatrix\right)^{-1}
  (\alpha_{1}-(x+iy\alpha_{2}))^{-1}b\right)^{-1}\\
&\quad\quad\quad\quad\quad\left(\Bbb I +a\alpha_{2}^{-1}
  \left(\alpha_{1}-(x+iy\alpha_{2})\right)^{-1}b\right);
\endaligned \tag{4.4}$$
when $\lambda=1$, $S=E_{1}$ is the sought expression of Theorem A.

\remark{Remark {4.5}} We can verify directly that these solutions are unitary
using the reality conditions.  {}From reality,
$$\align S^{*} & =\one
+b^{*}(\alpha_{1}^{*}-(x-iy\alpha_{2}^{*}))^{-1}\alpha_{2}^{*-1}a^{*}\\
& =\one-a \pmatrix & \one\\ \one & \endpmatrix
\left(\pmatrix & \one\\ \one & \endpmatrix\alpha_{1}
\pmatrix & \one\\ \one & \endpmatrix\right.\\
& \qquad \left.-\left(x+iy\pmatrix & \one\\ \one & \endpmatrix \alpha_{2}
\pmatrix & \one\\ \one & \endpmatrix\right)\right)^{-1}
\pmatrix & \one\\ \one & \endpmatrix \alpha_{2}^{-1}
\pmatrix & \one\\ \one & \endpmatrix
\pmatrix & \one\\ \one & \endpmatrix b\\
& =\one - a(\alpha_{1}-(x+iy\alpha_{2}))^{-1}\alpha_{2}^{-1}b.\endalign$$
Let
$$\align X_{1}&=(\alpha_{1}\alpha_{2}-x\alpha_{2}-iy\alpha_{2}{}^{2})\\
X_{2}&=(\alpha_{2}\alpha_{1}-x\alpha_{2}-iy\alpha_{2}{}^{2}),\endalign$$
the monad equation $[\alpha_{1},\alpha_{2}]+ba=0$ is equivalent to
$X_{1}-X_{2}=-ba$.
Using this notation,
$$\align S^{*}S & =\one +a(X_{1}{}^{-1}-X_{2}{}^{-1})b
-aX_{2}{}^{-1}baX_{1}{}^{-1}b\\
& = \one +a(X_{1}{}^{-1}-X_{2}{}^{-1})b
-aX_{2}{}^{-1}(X_{2} -X_{1} )X_{1}{}^{-1}b\\
& =\one.\endalign$$
\endremark

\remark{Remark 4.6} We can simplify the uniton equation as well.
But we can only show
that it is satisfied when $\gamma=0$.  Now let
$X=\alpha_{1}-(x+iy\alpha_{2})$, so
$$\gather S=\one+a\alpha_{2}^{-1}X^{-1}b\\
S^{-1}=S^{*}=\one-aX^{-1}\alpha_{2}^{-1}b.\endgather$$
We calculate
$$\gather \dd{x}S=a\alpha_{2}^{-1}X^{-2}b\\
S^{*}\dd{x}S=aX^{-1}\alpha_{2}^{-1}X^{-1}b\\
\dd{x}(S^{*}\dd{x}S)=aX^{-2}\alpha_{2}^{-1}X^{-1}b+aX^{-1}
\alpha_{2}^{-1}X^{-2}b\\
\dd{y}S=ia\alpha_{2}^{-1}X^{-1}\alpha_{2}X^{-1}b\\
S^{*}\dd{y}S=iaX^{-2}b\\
\dd{y}(S^{*}\dd{y}S)=-aX^{-1}\alpha_{2}X^{-2}b-aX^{-2}\alpha_{2} X^{-1}b.
\endgather$$
So
$$\text{uniton equation}=aX^{-2}\left\{\alpha_{2}^{-1}X
+X\alpha_{2}^{-1}-\alpha_{2}X-X\alpha_{2}\right\}X^{-2}b.$$
When $\gamma=0$, so that $\alpha_{2}$ is diagonal, the middle factor
is zero.  But it is not zero in general (use the example of a
$\operatorname{U}(3)$ monad in section 5).
\endremark

\head 5. Example of a $U(3)$ uniton\endhead

We have two ways of describing
holomorphic vector bundles, in terms of clutching functions
or in terms of monad data.  The monad approach should be
thought of as an attempt to reduce to the case we understand, line
bundles, which form the commutative part of the category of bundles.  Our
tool is the resolution and the homological algebra surrounding it.  As we
have seen, it is an effective tool, but not without its drawbacks.
  On a practical level, we know that time invariance is
effectively a restriction on the bundle and not an extra structure.  This
was clear and testable in terms of the clutching data in the example in
\cite{An1}, but not obvious for monads.

Ideally, one would construct a dictionary to translate from monad to
clutching data and back and take advantage of both representations.  This
was one motivation for trying to directly construct the monad data
corresponding to the clutching data in \cite{An1}.  We were able
to do this for this example, but it's not clear how to do it in
general.  The procedure is as follows.

\smallpagebreak
Beilinson's theorem on the existence of a monad representative for ${\Cal
V} \mapsto X$ is proven by defining a resolution of $\pi_{1}^{*} {\Cal V}
\otimes{\Cal O}_{\Delta}$ on $X\times X$ whose cohomology is the bundle to be
represented.  We can push the resolution down to $X$ via two projections.
The first push-down gives back $\Cal{V}$.
In good cases, the spectral sequence for the other direct
image doesn't degenerate to the bundle, but has a nontrivial $E_{2}$ term
whose entries are twisted cohomology groups of $\Cal{V}$ which form the monad
(see \cite{OSS}).  It's not surprising that this can be done for
generically trivial bundles on rational surfaces, since such
bundles can be pulled back to $\Bbb{P}^{2}$ via generic birational
equivalences.  The author's thesis contains the details of this for
$\T$.  We omit the justification as we only need it for this example.

To make the computations easier, cut out a neighbourhood of
$\{\lambda=0\}$ in $\T$ and glue it into the trivial bundle on
$(x,z)\in\Bbb{P}^{1}\times\Bbb{P}^{1}$.  Let ${\Cal Z}(p,q)={\Cal Z}
\otimes{\Cal O}(\pi_{1}^{-1}(p\text{ points})
+\pi_{2}^{-1}(q\text{ points}))$.
The Leray spectral sequence
for $H^{*}(\Bbb P^1\times\Bbb P^1,\Cal V(p,q))$ gives  `natural' bases
for \v{C}ech cohomology.  Computing the direct image of the maps in the
resolution (i.e.\ their natural transforms) means writing $H^1(x)$, etc,
in terms of the chosen bases.  The direct image of $1\otimes
x_{2}-x_{1}\otimes 1$ is then $(H^1(x)-xH^1(1))$, where tensor product and
addition of maps become matrix multiplication and addition.  The
resulting $\Bbb P^{1}\times\Bbb P^{1}$ monad (for the   jump
at $P_0$ only) is
$$\CD
\Cal O(-1,-2)^{4} @>{\left(\smallmatrix x&0&0&0\\ 0&x&0&0\\ -1&0&x&0\\
0&-1&0&x
\endsmallmatrix\right)} >{\left(\smallmatrix 0&0&0&0\\ 0&-z&0&0\\
z&0&0&0\\ z&-1&0&0\\
0&0&0&-z\\ 0&0&z/2&1/2\\ 0&0&-1&0\endsmallmatrix\right)} >
{\matrix \Cal O(-1,-1)^{4}\\ \oplus\\ \Cal O(0,2)^{7}\endmatrix} @>
{\left(\smallmatrix
0&z&0&0\\ 0&0&0&z\\ 0&0&-z/2&-1/2\\ -z&1&0&0\endsmallmatrix\right)}>
{\left(\smallmatrix
-1&x&0&0&0&0&0\\ 0&-1&0&0&x&0&0\\ 0&0&-1&1/2&0&x&0\\
0&0&0&x&0&0&0\endsmallmatrix\right)}> \Cal O(0,-1)^{4}.
\endCD$$
We want a monad on $\Bbb{P}^2$, so make the substitutions
$$\align x&=X/Y \\ z&=WY/(X+Y)^{2}\endalign$$
coming from the birational equivalence in figure 1.
(This is well-defined on the blown-up variety at the top of {figure 1}.)
Before we blow down, we need to make a monad transformation given by
$$(1+\lambda)^{2}\Bbb I_{4}\qquad \pmatrix (1+\lambda)^{2}\Bbb I_{4}&\\
&\Bbb I_{7}\endpmatrix\qquad \Bbb I_{4}$$
(which again makes sense on the double blow up).  We then use the group of
isomorphisms of the $\Bbb{P}^{2}$ monad to put the resulting monad into normal
form yielding normalised data
$$\matrix \gamma=\pmatrix 0&0&1&0\\ 0&0&0&1\\ 0&0&0&0\\ 0&0&0&0
\endpmatrix & \alpha'_{1}=\pmatrix 0&-2&0&4\\ 0&0&0&0\\ 0&0&0&2\\
0&0&0&0
\endpmatrix\\
a'=\pmatrix 0&0&0&0\\ 0&0&0&2\\ 2&0&0&0 \endpmatrix & b'=
\pmatrix 0&4&0\\ 4&0&0\\ 0&0&0\\ 2&0&0 \endpmatrix \endmatrix.$$
Plugging this into the formula in Theorem A, we get the same solution we
obtained in \cite{An1}.

\head {6.} Interpretation of Uniton Number \endhead

We will now show that uniton number corresponds to the length of the polar
jumping line.  Recall that $S\in \operatorname{Harm}(\Bbb{S}^{2},
\operatorname{U}(N))$ has uniton number $n$ if $S$ admits an extended
solution of the form
$$\widetilde{E}_{\lambda}=T_{0}+\lambda T_{1}+\dots \lambda^{n}T_{n}$$
and this is the shortest possible such solution.  We will assume, without
loss of generality, that $\widetilde{E}_{\lambda}$ is in Uhlenbeck normal form,
i.e.\
$$\operatorname{span}\{\operatorname{im}T_{0}(z):z\in\Bbb{C}\}=\Bbb{C}^{N}.
$$
For our purposes, the best definition of the length of
$\Cal{V}|_{P_{0}}$ is as the largest $l$ such that there exists a map
$$\Cal{O}_{P_{0}^{(l)}}(1)@>s>> \Cal{V}_{P_{0}^{(l)}},\ s|_{P_{0}}\ne0,$$
equivalently
$$\Cal{V}_{P_{0}^{(l)}}@>s>>\Cal{O}_{P_{0}^{(l)}}(-1),\ s|_{P_{0}}\ne0.$$
In terms of our monad, this is the smallest $l$ such that $\gamma^{l+1}=0$.

\demo{Proof} The first characterisation follows from the observation that if
$p_{1}$ and
$p_{2}$ are two points on $P_{0}\subset\Bbb{P}^{2}$, then
$$\operatorname{ker}K_{p_{2}}\circ J_{p_{1}}= \{s\in
H^{0}(P_{0},\Cal{V}): s(p_{1})=0\}\cong H^{0}(P_{0},\Cal{V}(-1)),$$
and that this works just as well for $P_{0}^{(l)}$ in which case all objects
 are defined
over $\Bbb{C}[\lambda]/(\lambda^{n+1})$ instead of over $\Bbb C$.  Since
the dual bundle is
described by the transposed monad, this gives the equivalence.  See
\cite{Ti} for another definition of length. \enddemo

The bundle $\Cal{V}$ comes with two sorts of framings.  The first is
a framing of $\Cal{V}$ over $G_{\infty}$ and
$\widetilde{T\Bbb{C}}{}^{*}$ (the union of the nonpolar fibres), call it
$g$.  The other is a family of framings $f_{y}$ over $G_{y}$ where $y$ is
in a neighbourhood of the real sections in $\Bbb{C}^{3}$.  They are
related by the change of frame $E_{\lambda}(y)$
$$g=E_{\lambda}(y)\cdot f_{y}$$
(which is only defined over $\{\lambda\in\Bbb{C}^{*}\}\cap G_{y}$).

Among extended solutions, $E_{\lambda}$ is determined by the property
$E_{\lambda}(\infty)=\Bbb{I}$.  This implies
$$E_{\lambda}(y)=\widetilde{E}_{\lambda}(\infty)^{-
1}\widetilde{E}_{\lambda}(y).$$
In particular, the sections
$$\widetilde{E}_{\lambda}(y)\cdot
f_{y}=\widetilde{E}_{\lambda}(\infty)\cdot  g$$
are holomorphic on a neighbourhood of $P_{0}$ in
$\T$ and have full rank away from $P_{0}$.

Since $E_\lambda$ is an analytic function on $\Bbb{S}^2$, it extends to a
holomorphic function on an open subset of $\Bbb{C}^2$, which we extend
trivially to $\Bbb{C}^3$.  This is what we mean by $E_{\lambda}(y)$.
By $G_z$ we will mean $G_{(z,0,\bar{z})}$.

Since
$\widetilde{E}_{\lambda}$ is in Uhlenbeck normal form, we can find
$z_{0}$, $z_{1}$ such that
$$T_{n}(z_{0})^{*}T_{0}(z_{1})\ne0.$$
The reality condition on $E_{\lambda}$
$$E_{\lambda}^{-1}=(E_{\bar{\lambda}^{-1}})^{*}=\lambda^{-n}T_n^* +\dots
+\lambda^{-1}T_1^*+T_0^*$$
implies that the $N$ sections
$$\lambda^{n}\widetilde{E}_{\lambda}(z_{0})^{-1}\widetilde{E}_{\lambda}(y)
\cdot
f_{y}=\lambda^{n}\widetilde{E}_{\lambda}(z_{0})^{-
1}\widetilde{E}_{\lambda}(\infty)
\cdot  g$$
have zeros along $P_{0}^{(n-1)}\cap G_{z_{0}}$, but are not all zero on
$P_{0}\cap G_{z_{1}}$.  This shows that $l\geq n-1$ (by the first
characterisation of length).

Assume, on the other hand, that $l> n-1$ and let $s$ be a section
$$\Cal{V}@>s>>\Cal{O}_{P_{0}^{(l)}}(-1);\ \  s|_{P_{0}}\ne0.$$
We can find a $y_{0}\in\Bbb R^3$ and a $v\in\Cal{V}_{P_{0}\cap G_{y_{0}}}$
such
that $s(v)\ne0$.  It follows that the sections
$$\lambda^{n}E_{\lambda}(y_{0})^{-1}E_{\lambda}(y)\cdot  f_{y}$$
on $P_{0}^{(n)}\cap G_{y_{0}}$ are not all mapped to zero by $s$, but
since their image is a section of $\Cal{O}_{P_{0}^{(n)}}(-1)\cong\Cal O_{P_0}
(-1)^{\oplus(n+1)}$, they
must be zero.  We conclude that $n=l+1$.

\head 7. Real triviality  \endhead

A bundle over $\Bbb P^1$ with zero first Chern class is either trivial or
contains a positive subbundle.  To show $\Cal{V}$
$|_{G_{\text{real}}}$ is trivial
we will show $H^0(G_{\text{real}},\Cal V)$ contains no sections with isolated
zeros.

Recall that $\Cal V$ is trivial on all real sections iff it is trivial on
the $\sigma$-invariant lines $\{\bar{z}X-zY+2W=0\}$.  Let $\overline{pq}$
be such a line, with $\sigma(p)\ne q$.

We need the
\proclaim {Lemma {7.1}} There exist hermitian metrics on $\Bbb C^k$ and $\Bbb
C^{2k+N}$
such that
$$(J_{\sigma(p)})^*=K_p.$$
\endproclaim
This allows us to calculate
$$\align
H^0(\overline{pq},\Cal V)& \cong\ker K_p \cap\ker K_q\\
&=(\operatorname{im} K_p^*)^\perp \cap\ker K_q\\
&=(\operatorname{im} J_{\sigma(p)})^\perp \cap\ker K_q\endalign$$
which implies
$$\align
\{s\in H^0(\overline{pq},\Cal V): s(\sigma(p))=0\}&
\cong\ker K_p\cap\ker K_q \cap \operatorname{im} J_{\sigma(p)}\\
&=\{0\}.\endalign$$
It follows that $\Cal V|_{\overline{pq}}$ is trivial.

\demo {Proof of Lemma}
Assume $J$ and $K$ are in the real form specified in Theorem B.  Let
$$\Omega\deq\pmatrix &\pmatrix&\one\\ \one \endpmatrix \\
\pmatrix&\one\\ \one \endpmatrix \\
& & \one
\endpmatrix\quad\text{and}\quad\omega\deq\pmatrix&\one\\ \one \endpmatrix,$$
and define hermitian forms on $\Bbb C^{2k+N}$ and $\Bbb C^k$ by
$$<v_1,v_2>\deq\overline{v_2}^t \Omega v_1\ \text{ and }\
<u_1,u_2>\deq\overline{u_2}^t \omega u_1.$$
Since
$$\matrix \pmatrix&\one\\ \one \endpmatrix\alpha_2
  \pmatrix&\one\\ \one \endpmatrix=-\overline{\alpha_2}{}^t &
\pmatrix&\one\\ \one \endpmatrix b=\bar{a}^t\\
\pmatrix&\one\\ \one \endpmatrix\alpha_1
  \pmatrix&\one\\ \one \endpmatrix=\overline{\alpha_1}{}^t &
a\pmatrix&\one\\ \one \endpmatrix=\bar{b}^t\endmatrix$$
$\overline{\sigma(X+Y)}=X+Y$, $\overline{\sigma(X-Y)}=-(X-Y)$
and $\overline{\sigma(W)}=-W$, we have
$$K_p=\overline{\Omega J_{\sigma(p)}\omega}^t=(J_{\sigma(p)})^*$$
as required.\enddemo

\head {8.} Time translation\endhead

Because \roster
\item $\Cal V|_{\text{nonpolar fibres}}$ is holomorphically trivial,
\item $\delta_{t}:\T\to\T$ preserves the fibres,
\item $\delta_{t}$ fixes $G_\infty$ and its lift $\tilde{\delta}_t$ fixes
$\Cal V|_{G_\infty}$, and
\item the union of the nonpolar fibres is a Zariski open set,
\endroster
the lift of time translation to the bundle is unique if it exists.
It therefore makes sense to refer to
`time-invariant' holomorphic bundles.  We must identify these
bundles.

Since the real structure, $\tilde{\sigma}$, mirrors the holomorphic
structure at $P_0$ and $P_\infty$, the same argument shows that time
invariance is a
local property, depending on the framed holomorphic structure in a
neighbourhood
of either pole.  In terms of monads, we may work with half the
monad (as in \S{3.7}) and look for conditions under which time translation
exists as
a
monad isomorphism.

In section 3 we fixed a birational equivalence of $\T$ and $\Bbb{P}^{2}$
which was convenient for deriving the formula.  In fact, any
equivalence by which (complete) uniton bundles could be pushed forward
must send
time translation to a discontinuous group of transformations on
$\Bbb{P}^{2}$.  Working with half the monad (and half the bundle,
$\Cal V_{\text{north}}$ of \S{3.7}), we may blow up $\hat{\eta}=0,
\lambda=\infty$, instead.  The resulting equivalence,
$$Y\lambda=X,\qquad Y\eta=W,$$
sends $\delta_{t}:\lambda\mapsto \lambda, \eta\mapsto\eta+\lambda t$ to
a {\it linear} transformation of $\Bbb{P}^{2}$:
$$X\mapsto X,\quad Y\mapsto Y \text{ and } W\mapsto W+Yt.\tag{8.1}$$

Since time invariance is a local property, it
doesn't matter which birational equivalence we choose as we get equivalent
monads.  So we may as well assume we are working with the $1/2$ monad
$(\alpha'_{1},\alpha'_{2},a',b')$.  We are looking for a map $\Cal{V}\to
\delta_t^*\Cal{V}$.  Such a map exists iff the respective monads are
equivalent under the group action.  Pulling $J$ and $K$ back by
$\delta_{t}$ \thetag{8.1} disturbs the normalisation we have chosen.
Multiplying $\Cal{O}^{\oplus(2k+N)}$ by
$$\pmatrix \one&t\one&\\&\one&\\&&\one
\endpmatrix$$
restores the normalisation, but sends
$$\align \alpha'_{1}&\mapsto \alpha'_{1}-2t\gamma\\
\alpha'_{2}&\mapsto \alpha'_{2}\\
a'&\mapsto a'\\
b'&\mapsto b'.\endalign$$
The two are equivalent iff there exists a one parameter subgroup
$G(t)\in\operatorname{Gl}(k/2)$ which fixes $\gamma$, $a'$ and $b'$, and
sends $\alpha'_1$ to $\alpha'_1-2t\gamma$.  Infinitesimally, this says there
exists $g\in\operatorname{gl}(k/2)$ such that
$$[g,\alpha'_{1}]=\gamma,\quad [g,\gamma]=0,\quad gb'=0,\quad
a'g=0.\tag{8.2}$$
As we have said, $g$ is unique up to the group action.  When it exists,
 it is the
$\delta$ of Theorem B.

\proclaim{Lemma {8.3}} No monad with
$$\gamma\in \left\{
\pmatrix N_n
\endpmatrix,
\pmatrix N_2 \\ & 0
\endpmatrix,
\pmatrix N_m \\ & Z_r
\endpmatrix
: n>1,m>2,r>0\right\}$$where
$$N_n\deq\pmatrix 0 & 1 \\
           & 0 & \ddots\\
           &    & \ddots & 1\\
           &    &        & 0
\endpmatrix\in\operatorname{gl}(n),\text{ and }
Z_r\deq\pmatrix0\endpmatrix\in\operatorname{gl}(r)$$
exists.
\endproclaim
\demo{Proof} Let $\alpha=\alpha_1'$, $\gamma$, $a=a'$, $b=b'$ be monad
data for $\Cal{V}_{\text{north}}$, and let $\delta$ and $g$ represent the
infinitesimal time translation and group action respectively.
Assume $\gamma$ of the form above.  We
 first show that $n\leq 2$.

If $[g,\gamma]=0$, then
$$g=\pmatrix g_{11} & g_{12}\\ g_{21}& g_{22}
\endpmatrix$$
where $g_{11}\in\operatorname{gl}(n)$ is upper-diagonal Toeplitz,
$g_{12}$ is zero except in its first row and $g_{21}$ is zero except in
its last column.  The same applies to $\delta$.

\smallpagebreak\noindent{\it Case I.}
If $a_1\ne0$ (the first column), we can assume
(use the group action to insure) that $(a_1,a_j)=0$, for all $j>1$.  With this
normalisation $ag=0$ implies $g_{11}=0$, $g_{12}=0$.  These restrictions
apply as well to $\delta$.  The condition $[\alpha,\delta]=\gamma$ now implies
$n\le2$, and when $n=2$, that the last
column of $\delta_{21}$ is not zero, while the second row of $\alpha$
is zero.  Nondegeneracy then implies that the second row of $b$ is
not zero, but $\delta b=0$ then implies $\delta=0$, which contradicts the
assumption that $[\alpha,\delta]=\gamma$.

\smallpagebreak\noindent{\it Case II.}
Assume that the first $n_a\ge1$ columns of $a $
are zero, as are the first $n_g$ columns of $g$.
Then $ag=0$ implies $n_a+n_g\ge n$.

If $n_a\ge2$, $[\alpha,\gamma]=ba$ implies that the first column of
$\alpha$ must be zero, violating nondegeneracy.  If, on the other hand,
$n_g\ge n-1$,
$[\alpha,\delta]=\gamma$ is only possible when $n\le 2$.

We are left with the case $n=2$, $a_1=0$, and $b_2=0$ (by symmetry).  The
monad equation and nondegeneracy imply
$$\align\alpha&=\pmatrix\alpha_{11} &\alpha_{12}&\alpha_{13}\\
0 &\alpha_{22}&\alpha_{23}\ne0\\
\alpha_{31}\ne0& \alpha_{32}&\alpha_{33}
\endpmatrix,\quad
b=\pmatrix b_1\ne0\\ 0\\b_3\ne0
\endpmatrix,\\
a&=\pmatrix 0 &a_2\ne0 &a_3\ne0
\endpmatrix,\endalign$$
and $a\delta=0$, $\delta b=0$ imply
$$\delta = \pmatrix0 & \delta_{12} & 0\\
0 & 0 & 0\\
0 & 0 & 0
\endpmatrix;$$
$([\alpha,\delta])_{32}=\gamma_{32}=0$ implies $\delta=0$ which
contradicts $[\alpha,\delta]=\gamma$.\hfill\qed
\enddemo
It is possible to extend this argument,
but not (it seems) without degenerating into a comedy of cases.  This
much implies
\proclaim{Corollary {8.4}}Two-unitons have normalised energy at least
four.  This bound is sharp (see \S5).\endproclaim

\head 9. Monads with $\gamma=0$\endhead

In the $\gamma=0$ case, we have just seen that time invariance poses no
additional constraints on the monad.  We will now use the residual action
of $G\in\operatorname{Gl}(k/2)$ and the nondegeneracy condition to put
$\alpha'_{1}$, ${a'}$ and ${b'}$ into a normal form.
In the next section we will then specialise to the well-known
$\operatorname{U}(2)$ unitons
and show how the monad data is equivalent to the pole and principal part
description of based rational maps.

Again, we choose to work with half the monad.
Since $\operatorname{Gl}(k/2)$ acts on $\alpha'_{1}$ by conjugation,
we can put
$\alpha'_{1}$ into
Jordan
normal form.
If we agree to a lexicographical ordering of $\Bbb{C}$, we can fix the
order of the Jordan blocks up to a permutation of blocks with the same
rank and
eigenvalue.
Unless $\alpha'_{1}$ is diagonalisable with distinct eigenvalues, its Jordan
form
has a nonzero stabiliser which continues to act on ${a'}$ and
${b'}$.
Once we divide $\alpha'_{1}$ into groups of Jordan blocks, each group with a
distinct eigenvalue and rank, we get a corresponding decomposition of
$\Bbb{C}^{k/2}$ into $\alpha'_{1}$-invariant subspaces.
They are the invariant subspaces of the stabiliser (under the standard
action of $\operatorname{Gl}(k/2)$ on $\Bbb{C}^{k/2}$), so we may
consider them one at a
time.

Assume for the moment that $\alpha'_{1}$ has exactly $m$ Jordan blocks of size
$l$.  It is convenient to assume
$$\alpha'_{1}=\pmatrix
\nu&&1&&\\&\ddots&&\ddots&\\&&\ddots&&1\\&&&\ddots&\\
&&&&\nu\endpmatrix,$$
with $lm$ $\nu$'s and $m(l-1)$ $1$'s, because the
stabiliser is then in block form:
$$\operatorname{Stab}_{\alpha'_{1}}\operatorname{Gl}(k/2)=\left\{ \pmatrix
g_{0}&g_{1}&\dots&g_{l-1}\\&\ddots&\ddots&\vdots\\&&\ddots&g_{1}\\
&&&g_{0}\endpmatrix :
\aligned g_{0}&\in\operatorname{Gl}(m)\\
g_{i}&\in\operatorname{gl}(m)\\
i&=1,\dots ,l-1\endaligned\right\}.$$

The injectivity of $J$ imposes an independence condition on
${a'}_{}$.  At $\lambda=0,\eta=-\nu$, $\one\eta+\alpha'_{1}$ is singular;
the first $l$ columns are zero, so if $u\in\Bbb C^{m}$
$$
\pmatrix u\\0\endpmatrix \in\pmatrix \Bbb{C}^{m}\\\oplus\\\Bbb{C}^{m(l-1)}
\endpmatrix
\matrix J|_{(\lambda,\eta) =(0,-\nu)}\\
\mapstochar\relbar\mathrel{\mkern-4mu}\relbar
\mathrel{\mkern-4mu}\relbar\mathrel{\mkern-4mu}\relbar\mathrel{\mkern-
4mu}\relbar
\mathrel{\mkern-4mu}\relbar\mathrel{\mkern-4mu}\relbar\mathrel{\mkern-
4mu}\relbar
\mathrel{\mkern-4mu}\relbar\mathrel{\mkern-4mu}\rightarrow\endmatrix
  \pmatrix (0)\\ {a'} \pmatrix u\\0\endpmatrix
           \endpmatrix
    \in \matrix\Bbb{C}^{2ml}\\ \oplus\\\Bbb C^{N}\endmatrix.$$
Since $J$ is injective, the first $m$ columns of $a'$ must be independent.
Similarly, the surjectivity of $K$ at the same point implies that
the last $m$ rows of $b'$ are independent.

Since $g_{0}$ acts on $({a'})_{1},\dots({a'})_{m}$ by right
multiplication, ${a'}$ (and ${b'}$) define points on a
Grassmannian:
$$\left\{m \text{ independent } N\text{-vectors }\right\}/\operatorname{Gl}
(m)\cong\operatorname{Gr}_{m}(\Bbb C^{N}).$$
If $\alpha'_{1}$ has Jordan blocks of different rank but the same eigenvalue,
we get an element of a Stiefel manifold, a space
of linearly independent subspaces with prescribed ranks.

In any case, within each invariant block, we can use the action of
$g_0$ to put the last $m$ rows of ${b'}$ into
some normal form, say the one which gives coordinates on the Schubert
cycles \cite{{GH},1.5}.
This reduces the stabiliser to the subgroup
$\{(\one, g_{1},\dots g_{l-1})\}$.

There are two `standard' metrics on $\Bbb{C}^N$ each of which can be used
to
put ${b'}$ into a normal form.  If we put the hermitian metric on
$\Bbb{C}^N$,
we can use the remaining action of the stabiliser to put the first
$m(l-1)$ columns of $b'$ perpendicular to the rest
$$\left\{b^{\prime }_{1},b^{\prime }_{2},\dots,b^{\prime }_{m(l-
1)}\right\}
\subset\left\{b^{\prime }_{m(l-1)+1},b^{\prime }_{m(l-1)+2},\dots,
b^{\prime }_{ml}\right\}^{\perp_{\text{hermitian}}}.$$
This shows that the data ${b'}$ describe a point in the orthogonal
bundle to the universal bundle contained in
$\underline{\Bbb{C}}^N_{Gr_{m}(\Bbb C^{N})}$.

Alternatively, we can put the holomorphic Euclidean metric on
$\Bbb{C}^N$.  In this case we have to worry about null vectors, so the
same
procedure doesn't work.  Instead we have to consider the usual coordinate
patches of $\operatorname{Gr}_{A}({\Bbb{C}}^{N})$.  If the last rows of
${b'}$ are
in Schubert
cycle form, i.e.\
$$\matrix\matrix\phantom{0}&\phantom{\dots} &\phantom{0}&\pretend
i_1\haswidth0
      &\phantom{0}& \phantom{\dots} &\phantom{0}&\pretend i_2\haswidth
0&\phantom{0}& \phantom{\dots}
      &\phantom{0}&\pretend i_3\haswidth0& \phantom{0}&\dots\\
   &  & &\downarrow& &     & &\downarrow& &     & &\downarrow&
   &\endmatrix\\
  \pmatrix 0&\dots&0&1&*&\dots&*&0&*&\dots&*&0&*&\dots\\
     0&\dots&0&0&0&\dots&0&1&*&\dots&*&0&*&\dots\\
     0&\dots&0&0&0&\dots&0&0&0&\dots&0&1&*&\dots\\
     \vdots &\ddots& \vdots&\vdots &
     \vdots&\ddots&\vdots&\vdots&\vdots&\ddots&\vdots&\vdots&\vdots&\ddots
     \endpmatrix,\endmatrix$$
then we can make the first $m(l-1)$ rows of ${b'}$ perpendicular to
$\{e_{i_1},e_{i_2},\dots,e_{i_m}\}$ (in either metric).

Apply this procedure to each invariant subspace of
$\operatorname{Stab}_{\alpha'_{1}}\operatorname{Gl}(k/2)$ and we
are left with a normalised monad which uniquely represents the bundle
$E$, i.e.\ its stabiliser in the group of monad isomorphisms is the trivial
subgroup.

\remark{Remark {9.1}} The reason for apparently two methods of
normalising is simple.  The moduli space of allowed $a'$'s given a choice of
$\alpha_1'$ is a bundle over the appropriate Grassmannian or Stiefel
manifold.  As complex bundles, they are trivial and
admit global frames.  Holomorphically, however, they are not
and we have to make different normalisations above the various Schubert
cells.\endremark

\remark{Remark {9.2}} To see directly
that $\gamma=0$ monads give $1$-unitons, we
can simplify the
formula in Theorem A to something which looks like the holomorphic map
composed with the Cartan embedding.  We could then imitate the extended
solution given in \cite{Uhl}.
\endremark

\head 10. Example: $\operatorname{Harm}(\Bbb
S^{2},\operatorname{U}(2))$\endhead

This section treats the simplest case:  $\operatorname{U}(2)$ unitons
given by rank-two
uniton bundles.  It demonstrates  the formula
 and points out an important distinction between maps into symmetric
spaces and maps into $\operatorname{U}(N)$.  Inside $\operatorname{U}(N)$,
holomorphic and antiholomorphic maps into symmetric spaces look the same.

It is well known that harmonic maps $\Bbb{S}^2\to\operatorname{U}(2)$
factor through spheres and hence are
closely linked to  rational maps.
We will give an ahistorical proof that based unitons correspond to
rational maps ($\Bbb{P}^1\to\Bbb{P}^1$), and show that the action of
$\operatorname{U}(2)$ on
$\operatorname{Harm}^{*}(\Bbb{S}^{2},\operatorname{U}(2))$
by conjugation
($S\mapsto USU^*$) corresponds to the usual
$\operatorname{Gl}(2)$ action on rational maps, i.e.\ the correspondence
is equivariant.
We then show that this is the same map as given by
Theorem A and we find that real triviality
\thetag{3.10} is implied by the other monad conditions in the
$\operatorname{Harm}(\Bbb{S}^{2},\operatorname{U}(2))$ case.  This
implies that $\operatorname{U}(2)$ monads always have $\gamma=0$.
Of course, this is not the way we would like to prove it.

\subhead 10.1 Rational maps\endsubhead
Since harmonic maps into $\operatorname{U}(2)$ are 1-unitons, they
have the form  $S=Q(\pi-\pi^\perp)$, for some
$Q\in\operatorname{U}(2) $, and
$\pi$, projection onto a holomorphic subbundle of
$\underline{\Bbb{C}^N}$.

We can see the decomposition $S=Q(\pi-\pi^\perp)$ as a composition of
three maps.
In the middle is the Cartan embedding
$$\Bbb{P}^1=\operatorname{Gr}_{1}(\Bbb C^{2})\overset I\to\hookrightarrow
\operatorname{U}(2)
:\pi\in\operatorname{Gr}_{1}(\Bbb C^{2})\mapsto\pi-\pi^\perp\in
\operatorname{U}(2).$$
In terms of $z\in\Bbb{P}^1=\operatorname{Gr}_{1}
(\Bbb{C}^{2})$,
$$\align\pi_z&=\frac1{1+z\bar{z}}\pmatrix 1\\z\endpmatrix\pmatrix 1\\
z\endpmatrix^*,\\
\pi_z^\perp&=\frac1{1+z\bar{z}}\pmatrix -\bar{z}\\1\endpmatrix\pmatrix
-\bar{z}\\
1\endpmatrix^*,\\ \intertext{so}
I(z)&=\pi_z-\pi_z^\perp=\frac1{1+z\bar{z}}\pmatrix1-z\bar{z}&2\bar{z}\\
2z&z\bar{z}-1\endpmatrix.\tag{10.2}
  \endalign$$
All holomorphic subbundles, $\pi\subset\underline{\Bbb{C}}^2$,
are
$\pi_f$ for some rational function $f:\Bbb{P}^1\to\Bbb{P}^1$.
Finally, we can left-translate harmonic maps $L_Q:S\mapsto QS$.
In other words, a general one uniton is a composition
$$\Bbb{P}^1@>f>>\Bbb{P}^1\overset
I\to\hookrightarrow\operatorname{U}(2)@>L_Q>>
\operatorname{U}(2) $$
where $f\in\operatorname{Rat}\Bbb{P}^1$ and $Q\in\operatorname{U}(2)$.
To what extent is it unique?

The group $\operatorname{U}(2)$ acts on itself by left translation, and
in turn on the totally-geodesically-embedded spheres.  The translates
of $I(\Bbb P^{1})$ are $\operatorname{U}(2)/\operatorname{Stab}
_{I(\Bbb P^{1})}\operatorname{U}(2)$.  Since $I(\Bbb P^{1})=\{S\in
\operatorname{U}(2):S^{2}=\one\}$, and $\one=(QI(0))^{2}=(QI(1))^{2}
=(QI(i))^{2}$ iff $Q=\{\pm\one\}$, $\operatorname{Stab}_{I(\Bbb P^{1})}
\operatorname{U}(2)=\{\pm\one\}$.  So left translation on embedded
spheres is not faithful. Since $(-\one)I(z)=I(-1/\bar{z})$ is
orientation reversing, it is
faithful on oriented spheres, and hence on
$\operatorname{Harm}(\Bbb{S}^{2},\operatorname{U}(2))$, i.e.\ the
decomposition $S=Q(\pi-\pi^{\perp}):(\operatorname{U}(2)\times
\operatorname{Rat}(\Bbb P^{1}))\to
\operatorname{Harm}(\Bbb{S}^{2},\operatorname{U}(2))$ is unique.

Rational maps can be written as $p(z)/q(z)$, with $(p,q)=1$.
Their topological degree is given by $\max\{\deg p,\deg q\}$ (always
positive, because holomorphic maps preserve orientation).
They contain the based maps
$$\operatorname{Rat}^*=\left\{f\in\operatorname{Rat}:f(\infty)=0\right\}=
\left\{\frac{p(z)}{q(z)}:(p,q)=1,\deg p<\deg
q\right\}.$$
Various groups act on $\operatorname{Rat}$ via the action of
$\operatorname{PGL} (2)$ given by
$$\pmatrix A&B\\C&D\endpmatrix:\frac pq\mapsto\frac{Ap+Bq}{Cp+Dq}.$$
This map preserves degree because $\operatorname{Gl}(2)$
is connected and degree
components are disjoint.
$\operatorname{Rat}$ is a $\operatorname{Rat}^*$ bundle over $\Bbb{P}^1$,
given by
$$\text{rat}:\operatorname{Rat}\to\Bbb{P}^1:f\mapsto f(\infty).$$
$\operatorname{PGL}(2)$ acts on $\Bbb{P}^1$ by
$$\pmatrix A&B\\C&D\endpmatrix:[x,y]\mapsto[Ax+By,Cx+Dy],$$
making $\operatorname{Rat}\to\Bbb{P}^1$ an equivariant bundle.
$$\{P\in\operatorname{PGL}(2):P(\operatorname{Rat}^*)
=\operatorname{Rat}^*\}=\left\{\pmatrix1&0\\C&D\endpmatrix\right\}.$$

Conjugation acts on $\operatorname{Harm}^{*}(\Bbb{S}^{2},\operatorname{U}(2))$
and the isomorphism $\operatorname{Harm}^{*}(\Bbb{S}^{2},\Bbb{U}(2))
\cong\operatorname{Rat}\Bbb{P}^1$ is equivariant.
\proclaim{Claim {10.3}} If $p/q\in\operatorname{Rat}\Bbb{P}^1$ and
$U=\pmatrix A&B\\C&D\endpmatrix\in\operatorname{U}(2)$
then $$UI(p/q)U^*=I\left(\frac{Cq+Dp}{Aq+Bp}\right).\tag{10.4}$$\endproclaim
\demo{Hint}
Write
$$I(p/q)=\frac{1}{p\bar p+q\bar q}\pmatrix q&-\bar p\\p&\bar q\endpmatrix
\pmatrix 1\\&-1\endpmatrix\pmatrix q&-\bar p\\p&\bar q\endpmatrix^*.\qed$$
 \enddemo

\subhead 10.5 Monads\endsubhead
We will now relate monad description to the rational maps description.
The second Chern class gives a stratification
$$\operatorname{Harm}^{*}(\Bbb{S}^{2},\operatorname{U}(2))
=\bigcup_{k/2\in\Bbb{N}}
\operatorname{Harm}^{*}_{k/2}(\Bbb{S}^{2},\operatorname{U}(2)).$$
(The quantity $k/2$ is also the jumping type of $\Cal V|_{P_0}$).
Degree gives a stratification
$$\operatorname{Rat}=\bigcup_{j}\operatorname{Rat}_{j}.$$
The map $\operatorname{Harm}(\Bbb{S}^{2},\operatorname{U}(2))/
\operatorname{U}(2)\to\operatorname{Rat}$ given by
the formula of Theorem A preserves this stratification.

Recall the Jordan block normalisation of a $\gamma=0$
monad $M$, given by the data $\alpha'_{1}, {a'}, {b'}$.
{}{}From the monad equation, we saw that ${b'}{a'}=0$, but
both ${b'}$ and ${a'}$ must be nonzero if $J$ and $K$ are to be
injective and surjective respectively.
It follows that the column space of ${a'}$ and row space of
${b'}$ are one dimensional, and that they are perpendicular to one another
with respect to the Euclidean metric on $\Bbb{C}^2$.
{}{}From the discussion of the normalisation, we see that the Jordan
blocks of $\alpha'_{1}$ have distinct eigenvalues,
and the monad will be given by
$$\align
\alpha'_{1}&=\pmatrix \operatorname{Jordan}(j_1,e_1)\\&\ddots
\\&&\operatorname{Jordan}(j_L,e_L)\endpmatrix,\\
{a'}&=\pmatrix1\\c_1\endpmatrix(1,0,\dots,0,1,0,\dots),\tag{10.6}\\
{b'}&=\tilde{b}
\pmatrix-c_1&1\endpmatrix\endalign$$
generically, where $\tilde{b}_{j_1}$,$\tilde{b}_{j_1+j_2},\dots$
are not zero, and the first, $j_1+1\text{st},\dots$ columns of $a$
are not zero.

When $c_1=0$,
	$$\frac{1}{2}{a'}(\alpha'_{1}-2z)^{-1}{b'}=\pmatrix
	0&f\\0&0\endpmatrix\tag{10.7}$$
where
$$ f(z)= - \sum_{i=1}^L\sum_{j=1}^{j_i}\frac{\tilde{b}_{j_1+\dots+j_i+j}
(z-e_i)^{j_i-j}}
  {(z-e_i)^{j_i}}=\frac pq\tag{10.8}$$
is a based rational map of degree $k/2$.  We see that
$$\align S&=\frac1{1+f\bar f}\pmatrix1-f\bar f&2f\\-2\bar f&1-f\bar
f\endpmatrix.\endalign$$
This is equivalent to \thetag{10.2} under a change of frame
(basing condition).

\remark{Remark {10.9}}
We see from this calculation that {\it {all}\/} $\operatorname{U}(2)$
unitons have $\gamma=0$, as expected.
\endremark

\remark{Remark {10.10}}
{}Putting this into the determinant form of real
triviality \thetag{3.10}, we can
$$\align\thetag{3.10}&=\prod_{i=1}^{L}|(e_i+z)|^{2j_i}
    \left(1+\left|\sum_{i=1}^{L}\sum_{j=1}^{j_i}
       (e_i+z)^{-j}\tilde{b}_{j_1+\dots+j_i+j}\right|^i\right)\\
  &=|q|^2+|p|^2\endalign$$
where $f=p/q$ as above, which is not zero since $(p,q)=1$ by construction.
\endremark

\remark{Remark {10.11}}
The interpretation of the determinant condition $\thetag{3.10}=|q|^{2}
+|p|^{2}$ tells us what happens (in the $\gamma=0$ case) when it fails.
The points where it fails are the common zeros of $p$ and $q$.
Nondegenerate monads do not give rise to such polynomials.
In the degenerate case, i.e.\ on part
of the boundary, the poles and zeros of $f$ coalesce and the holomorphic
map changes degree as the harmonic map undergoes bubbling off.\endremark

\remark{Remark {10.12}}
Having worked out the simple $\operatorname{U}(2)$ case in detail, we can
already give an interesting example (\cite{BuGu}).
  One-unitons are known to be holomorphic or antiholomorphic.
The simplest such solution is $I(z)$ in our notation, which is the map into
$\Bbb{P}^1$ associated to the tautological subbundle of
$\underline{\Bbb{C}}^2_{\Bbb{P}^1}$.  The orthogonal complement of this
bundle is represented by $-I(z)$.  The subbundle of
$\underline{\Bbb{C}}^4_{\Bbb{P}^1}$, which is the sum of
the tautological bundle
in the first $\underline{\Bbb{C}}^2$ and its complement in the second
$\underline{\Bbb{C}}^2$, is not holomorphic in
$\operatorname{Gr}_2(\Bbb{C}^4)$, so one might expect it to be a two-uniton.
As a map into $\operatorname{U}(4)$, however, it is equivalent, under left
multiplication by a scalar, to a holomorphic map into
$\operatorname{Gr}_2(\Bbb{C}^4)$.  In fact, the uniton bundle decomposes
into two copies of the simplest nontrivial bundle, and hence must be a
1-uniton.  This also tells us that the map has energy $2$.\endremark

\head 11. Grassmannian solutions \endhead
The involution $S\mapsto S^{-1}$ is an involution of $\operatorname{U}(N)$,
whose fixed set has $N-1$ components corresponding to the various Grassmann
manifolds of subspaces of $\Bbb{C}^{N}$.
This fixed set is the image of the Cartan embedding which takes a subspace to
the difference of projection onto the subspace and projection onto its
orthogonal complement.  Because this is a totally geodesic embedding, harmonic
maps into Grassmannians are in bijection with unitons whose image lies in the
fixed set of $S\mapsto S^{-1}$.
To be able to say something about these maps, we will restrict ourselves to
based maps.  The usual basing condition, $S(\infty)=\one$, does not make sense
for such maps.  We must choose different basepoints in each component,
$$S(\infty)=\pmatrix \one\\&-\one\endpmatrix$$
where the signature of the basepoint distinguishes the image Grassmannian.
We have described an inclusion
$$\operatorname{Harm}^*(\Bbb{S}^{2},\operatorname{Gr}_j(\Bbb{C}^N))
\hookrightarrow\operatorname{Harm}^*(\Bbb{S}^{2},\operatorname{U}(N)).$$
We would like to identify intrinsically which uniton bundles correspond to
based maps into Grassmannians.

Let $S\in\operatorname{Harm}^*(\Bbb{S}^{2},\operatorname{Gr}_j(\Bbb{C}^N))$
and let
$E_{\lambda}$  be an extended solution for $S$ and assume that
$E_{\lambda}(\infty)=\pmatrix \one&\\&-\lambda\one \endpmatrix$.
Uhlenbeck shows that
$\widetilde{E}_{\lambda}=E_{-\lambda}E_{1}^{-1}$
is an extended solution for $S^{-1}$.  If $S=S^{-1}$, then these extended
solutions are equivalent (after left multiplication by a function of
$\lambda$), and hence determine the same uniton bundle, up to framing.

Now define an involution by
$$\mu^{*}\lambda=-\lambda,\quad \mu^{*}\eta=\eta \quad
(\text{equivalently }\mu^{*}z=z),\tag{11.1}$$
and fix generators of $\pi_1(G_\infty\cup P_{-1}\cup G_z \cup P_{\lambda})$
above which we calculate the monodromies which determine $E_\lambda$:
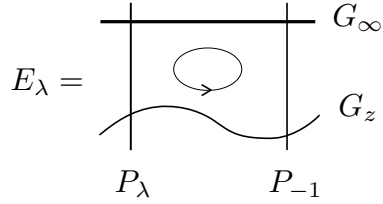
$$E_\lambda=\matrix\epsfxsize=70pt\epsffile{uniton1.epsf}\endmatrix
 \kern-68pt
 \matrix   \\ & & &  \kern-7pt  G_{\infty} \\
  \smallmatrix{}\  \\\ \\\  \\ \endsmallmatrix \\
   & &  & \kern-7pt G_z \\ \\
  P_\lambda&   \ \ \ \ \ &P_{-1}&  \endmatrix  $$
Then the formula for $\widetilde{E}_\lambda$ has the interpretation
$$\align E_{-\lambda}E_1^{-1}&=
  \matrix\epsfxsize=70pt\epsffile{uniton1.epsf}\endmatrix
 \kern-69pt
 \matrix \\  \\ \\ \\   \\
  P_{-\lambda}& \kern12pt  &P_{-1}&  \endmatrix
\circ\matrix\epsfxsize=70pt\epsffile{uniton3.epsf}\endmatrix
 \kern-69pt
 \matrix \\  \\ \\ \\   \\
  P_{1}& \ \ \ \ \   &P_{-1}&  \endmatrix
=\matrix\epsfxsize=100pt\epsffile{uniton2.epsf}\endmatrix
 \kern-95pt
 \matrix \\  \\ \\ \\   \\
  P_{\lambda} \kern26pt P_{-1} \kern19pt P_{1} \endmatrix\\
&=\matrix\epsfxsize=70pt\epsffile{uniton1.epsf}\endmatrix
 \kern-69pt
 \matrix \\  \\ \\ \\   \\
  P_{-\lambda}& \ \ \ \   &P_{1}&  \endmatrix
=\mu^*\pmatrix\epsfxsize=70pt\epsffile{uniton1.epsf}\endpmatrix
 \kern-78pt
 \matrix \\  \\ \\ \\   \\
  P_{\lambda}& \ \ \ \ \kern3pt   &P_{-1}&  \endmatrix
  =\mu^*E_\lambda.\endalign$$
Since $\widetilde{E}_\lambda$ determines the uniton bundle
and vice versa,
the uniton bundle for $S^{-1}$ is the $\mu$-pullback of the bundle for
$S$, up to the choice of framing.  To see the effect on the framing, note
that when $E_{\lambda}(\infty)= \pmatrix \one&\\&-\lambda\one
\endpmatrix$ (the difference between the chosen framing of
$\Cal{V}|_{G_{\infty}}$ and the canonical one) $\mu$ carries a frame
at $P_{-1}\cap G_{\infty}$ to $\pmatrix \one&\\&-\one \endpmatrix$ times
itself (after `transporting' it back using evaluation of
the canonical frame).

So uniton bundles ($\Cal{V},\Phi$) correspond to Grassmannian
solutions iff $\mu$ lifts to $\tilde{\mu}:\Cal{V}\to\Cal{V}$ and
the signature of $\mu^{*}\phi\phi^{-1}:\Bbb{C}^N\to\Bbb{C}^N$
determines the component (rank of the image Grassmannian).

\subhead {11.2} Monad characterisation \endsubhead
The bundle lift $\tilde{\mu}$ is a bundle map $\Cal{V}\to\mu^*\Cal{V}$,
but we want more than this, namely a map of monads.  On $\Bbb{P}^{2}$
monads are mapped to monads by linear maps.  The induced action of $\mu$
is linear if we take the birational equivalence from  \S{8},
i.e.\ the one in which we blow up and down on
$P_{\infty}$.  The induced action is
$$\matrix &X\mapsto -X\\ \mu:&Y\mapsto Y\\ &W\mapsto W.
\endmatrix$$
As in the case of time invariance, this means we must work with half the
bundle.  In what has gone before it has been convenient to frame our
bundles above $P_{-1}$.  To see that $\mu$-invariance is a local property
one should think of bundles framed over $G_{\infty}$.  The involution
$\tilde{\mu}$ may as well be thought of as an involution of
the space of frames of a given bundle.

If $(J,K)$ is a normalised monad representing $\Cal{V}_{\text{north}}$
$$\matrix J_{W}=\pmatrix \Bbb I\\0\\0\endpmatrix &
J_{X}-J_{Y}=\pmatrix 0\\ \Bbb I\\0\endpmatrix &
J_{X}+J_{Y}=\pmatrix \alpha_{1}\\ \alpha_{2} \\a\endpmatrix \\
\qquad&\qquad&\qquad\\
K_{W}=\pmatrix 0 & \Bbb I & 0 \endpmatrix &
K_{X}-K_{Y}=\pmatrix -\Bbb I & 0 & 0 \endpmatrix &
K_{X}+K_{Y}=\pmatrix -\alpha_{2} & \alpha_{1} & b \endpmatrix
\endmatrix\tag{11.3}$$
(primes omitted)  then the pulled-back monad
$$\gathered\mu^*J_{W}=\pmatrix \Bbb{I}\\0\\0
\endpmatrix\qquad \mu^*J_{X-Y}=\pmatrix -\alpha_{1}\\-\alpha_{2}\\-a
\endpmatrix\qquad \mu^*J_{X+Y}=\pmatrix 0\\-\Bbb{I}\\0
\endpmatrix\\
\mu^*K_{W}=\pmatrix 0&\Bbb{I}&0
\endpmatrix\qquad \mu^*K_{X-Y}=\pmatrix \alpha_{2}&-\alpha_{1}&-b
\endpmatrix\qquad \mu^*K_{X+Y}=\pmatrix\Bbb{I}&0&0
\endpmatrix\endgathered\tag{11.4}$$
represents $\mu^*\Cal{V}$.  If $G\in\operatorname{Gl}(k/2)$,
$F\in\operatorname{Gl}(N)$, then
$$  G\quad\quad\pmatrix G&-G\alpha_{2}^{-1}\alpha_{1}&-G\alpha_{2}^{-1}b\\
&-G\alpha_{2}^{-1}&\\&-Fa\alpha_{2}^{-1}&F
\endpmatrix\quad\quad G\alpha_{2}^{-1}
\tag{11.5}$$
(acting linearly on the three homogeneous bundles of the monad)
renormalises $\mu^*(J,K)$.  {}From this we deduce that
$\Cal{V}\cong \mu^*\Cal{V}$
iff there exist $G$ and $F$ such that
$$\matrix G\alpha_{2}^{-1}\alpha_{1}G^{-1}=\alpha_{1} &
Fa\alpha_{2}^{-1}G^{-1}=a\\ -G\alpha_{2}G^{-1}=\alpha_{2}&
-G\alpha_{2}^{-2}bF^{-1}=b,
\endmatrix\tag{11.6}$$
(primes omited) and $F$ represents the action on the framing;
in our case we require $F^{2}=\Bbb{I}$.

\remark{Remark 11.7} In the case of $1$-unitons, $\alpha_{2}$ is diagonal
and the monad equation reduces to $ba=0$.  Letting $F$ act on $\ker{b}$
by $-\Bbb{I}$ and on the complement by $\Bbb{I}$, and taking $G=\Bbb{I}$,
we see that all (based) $1$-unitons come from (based)
maps into Grassmannians.
\endremark

\subhead Thanks \endsubhead
I am greatful to Francis Burstall, Arleigh Crawford, Martin Guest,
Nigel Hitchin, Paul Norbury,
John Rawnsley, and especially Jacques Hurtubise
 for helpful discussions and friendly suggestions.

\enddocument